\documentclass[12pt]{amsart}
\usepackage{amsmath,amsthm}
\usepackage{amssymb}
\usepackage{euscript}           

\newtheorem{thm}{Theorem}[section]
\newtheorem{prop}[thm]{Proposition}
\newtheorem{cor}[thm]{Corollary}
\newtheorem{lemma}[thm]{Lemma}

\theoremstyle{definition}
\newtheorem{defn}[thm]{Definition}

\newtheorem{example}[thm]{Example}

\theoremstyle{remark}
\newtheorem{remark}[thm]{Remark}

\numberwithin{equation}{section}

\DeclareMathOperator{\colim}{colim\,}

\DeclareMathOperator{\ob}{ob}

\DeclareMathOperator{\Ch}{Ch}

\DeclareMathOperator{\Def}{Def}

\DeclareMathOperator{\mor}{mor}
\DeclareMathOperator{\Mor}{Mor}

\DeclareMathOperator{\Mod}{Mod}

\DeclareMathOperator{\sh}{Sh}

\DeclareMathOperator{\pre}{Pre}

\newcommand{\bd}{\textbf}

\newcommand{\ACC}{\mathsf{ACC}}

\newcommand{\TOPOI}{\mathsf{TOPOI}}

\newcommand{\HOMODEL}{\mathsf{HOMODEL}}

\newcommand{\cof}{\mathsf{cof}}

\newcommand{\inj}{\mathsf{inj}}

\newcommand{\dom}{\mathsf{dom}}

\newcommand{\ev}{\mathsf{ev}}

\newcommand{\bk}{\mathsf{bk}}

\newcommand{\codom}{\mathsf{codom}}

\newcommand{\Ex}{\mathsf{Ex}}

\newcommand{\exfty}{\mathsf{Ex}^\infty}
\newcommand{\Id}{\mathsf{Id}}
\newcommand{\cell}{\mathsf{cell}}

\newcommand{\Set}{\mathit{Set}}
\newcommand{\SSet}{\mathit{SSet}}

\newcommand{\Cat}{\mathit{Cat}}

\newcommand{\SRing}{\mathit{SRing}}
\newcommand{\Ring}{\mathit{Ring}}
\newcommand{\op}{{\text{op}}}
\renewcommand{\th}{{\text{th}}}
\newcommand{\mono}{\mathsf{mono}}
\newcommand{\isom}{\mathsf{iso}}
\newcommand{\AdRos}{Ad\'amek--Rosick\'y~\cite{adros}}

\newcommand{\topos}{\EuScript}
\newcommand{\cate}{\mathcal}

\newcommand{\A}{\cate A}

\newcommand{\B}{\cate B}

\newcommand{\C}{\cate C}
\newcommand{\tC}{\mathsf C}
\newcommand{\D}{\cate D}

\newcommand{\E}{\topos E}

\newcommand{\F}{\topos F}

\newcommand{\G}{\topos G}

\newcommand{\K}{\cate K}

\renewcommand{\L}{\topos L}
\newcommand{\M}{\cate M}

\newcommand{\cN}{\cate N}

\newcommand{\Q}{\topos Q}

\newcommand{\tS}{\mathsf S}

\newcommand{\Z}{\mathbb{Z}}
\newcommand{\tW}{\mathsf W}

\newcommand{\si}{{\Delta^{\text{op}}}}

\newcommand{\Xb}{X_\bullet}
\newcommand{\Yb}{Y_\bullet}
\newcommand{\lc}{\langle}                       
\newcommand{\rc}{\rangle}                       

\renewcommand{\preceq}{\preccurlyeq}		
\newcommand{\iso}{\cong}
\newcommand{\es}{\wedge}                        
\newcommand{\vagy}{\vee}                        
\newcommand{\tens}{\otimes}

\newcommand{\mon}{\rightarrowtail}
\newcommand{\epi}{\twoheadrightarrow}
\newcommand{\ra}{\rightarrow}
\newcommand{\lra}{\longrightarrow}

\newcommand{\inc}{\hookrightarrow}

\newcommand{\leri}{\leftrightarrows}

\newcommand{\three}[3]{\overset{#1}{\underset{#3}{#2}}}

\newcommand{\llra}[1]{\stackrel{#1}{\lra}}      
\newcommand{\xra}[1]{\xrightarrow{#1}}          

\usepackage[matrix,arrow,curve,cmtip]{xy}
\let\objectstyle=\displaystyle

\xymatrixcolsep{1.9pc}                     
\xymatrixrowsep{1.9pc}
\newdir{ >}{{}*!/-5pt/\dir{>}}             
\newdir{ |}{{}*!/-5pt/\dir{|}}             

\begin{document}

\title{Sheafifiable homotopy model categories}
\author{Tibor Beke}
\address{Department of Mathematics\\University of Michigan, Ann Arbor}
\email{tbeke@math.lsa.umich.edu}
\begin{abstract}
If a Quillen model category can be specified using a certain logical syntax
(intuitively, ``is algebraic/combinatorial enough''), so that it can be
defined in any category of sheaves, then the satisfaction of Quillen's
axioms over any site is a purely formal consequence of their being
satisfied over the category of sets.  Such data give rise to a functor from
the category of topoi and geometric morphisms to Quillen model categories
and Quillen adjunctions.
\end{abstract}
\maketitle

\section*{Introduction}
\noindent
The \emph{homotopy model category} of the title refers to the (``closed'')
model categories of Quillen~\cite{quil67}.  The intended meaning of
\emph{sheafifiable} is best illustrated by some examples:

\begin{example} (simplicial sheaves) \\   \label{jar}
Quillen's homotopy theory of simplicial sets can be extended to simplicial
sheaves over a site (i.e.\ category $\C$ with a Grothendieck topology $J$)
as follows.  Choose the cofibrations to be the monomorphisms.  Given a
choice of local basepoint for $X\in\sh(\C,J;\SSet)$, one can construct a
sheaf of homotopy groups over the basepoint.  See the details in
Jardine~\cite{jard87}.  A weak equivalence is a map that induces
isomorphisms on the homotopy sheaves for arbitrary local basepoints.  This
fixes the data for a Quillen model category, which is in fact simplicial
and proper.
\end{example}

\begin{example} (simplicial objects in a topos) \\   \label{joy}
A topos is a category equivalent to the category of sheaves on some site;
so a category of simplicial sheaves is a category of simplicial objects in
a topos.  In 1984, A.\ Joyal \cite{joyal} extended Quillen's homotopy
theory of simplicial sets to simplicial objects in a topos $\E$ as follows:
for $\Xb\in\E^\si$, its homotopy groups can be constructed --- by purely
categorical operations in $\E$ --- as objects (with algebraic structure)
over $X_0$. Let a weak equivalence be a morphism $\Xb\ra\Yb$ for which the
induced squares
\[ \xymatrix{ \pi_n(\Xb) \ar[d] \ar[r] & \pi_n(\Yb) \ar[d] \\
                  X_0 \ar[r] & Y_0 } \]
are pullbacks.  The rest is as above.
\end{example}

It is not hard to see that the two constructions prescribe identical
homotopy model structures for the same category.  For any topos $\E$, there
exist many (in fact, a proper class) of sites whose categories of sheaves
are all equivalent to $\E$; Joyal's result contains the additional
information that to all such sites, Ex.~\ref{jar} associates the same
homotopy theory.  (Historically, Ex.~\ref{joy} preceded Ex~\ref{jar}; a
detailed reworking of Joyal's proof appears in Jardine~\cite{jard96}.) 

\begin{example} (simplicial rings) \\   \label{siring}
Quillen~\cite{quil67} proves that the category of simplicial rings
possesses a homotopy model structure; the forgetful functor 
$\SRing\ra\SSet$ detects weak equivalences and fibrations.  He leaves
open the question whether the same is true for \emph{sheaves} of
simplicial rings.  The answer is affirmative: if $\Ring(\E)$ denotes the
category of ring objects of a topos $\E$, then with the help of the
canonical forgetful functor $\Ring(\E)\ra\E$ one can endow $\SRing(\E)$
with a model structure.  The remark in the paragraph following
Ex.~\ref{joy} applies here as well: since $\E^\si$ carries an intrinsic
homotopy structure (i.e.\ independent of having chosen a site for $\E$
first), the same is true for $\SRing(\E)$.
\end{example}

\begin{example} (small categories and categorical equivalences) \\
There is a ``folk model structure'' on $\Cat$, the category of small 
categories: the weak equivalences are the functors that induce an 
equivalence of categories, and the cofibrations the functors that are
injective on objects.  Joyal and Tierney~\cite{joytier90} proved that
(internal) category objects of any topos carry a homotopy model structure
with weak equivalences those functors that are (in the internal sense)
full, faithful and essentially surjective, and cofibrations those functors
that are monos on the object part.  (In the topos of sets, this specializes
to the case of $\Cat$ mentioned above.  It is distinct from the homotopy
theory of $\Cat$ introduced by Thomason~\cite{thom80} --- which, 
incidentally, sheafifies as well.)
\label{sstack}
\end{example}

Observe that in each case
\begin{itemize}
\item[(1)] One takes as granted a homotopy theory of structured sets, and
attempts to build one for sheaves of such structures.  The passage is not
arbitrary; for example, weak equivalences in the latter are locally weak
equivalences in the former.
\item[(2)] The homotopy structure is functorial in the topos, not merely 
in the site defining the topos.  (This is the difference between
Ex.~\ref{jar} and~\ref{joy}.)
\item[(3)] If $\E\llra{f}\F$ is a topos morphism\footnote{Throughout 
this article, topoi and their morphisms are understood in the sense of 
Grothendieck.  Every topos is a $U$-topos for a fixed universe $U$ whose 
elements are called simply \emph{sets}.}, then the inverse image functor
$f^*$ preserves the weak equivalences and cofibrations of the homotopy
model structure associated to $\F$, and induces a Quillen adjoint pair
between the corresponding model categories.  In particular, the ``stalk''
at a point $\Set\ra\F$ is the set-based homotopy model category one started
with; so, in this sense, one may think of the structure $\F$ is endowed
with as a ``sheaf of homotopy theories with constant stalk''.
\item[(4)] The homotopy model categories in question are all cofibrantly
generated.
\end{itemize}

\emph{Outline of this paper.} It is not hard to find a precise sense of
what was loosely referred to above as ``structured sets'' and ``sheaves of
such structures'', and to guess what additional criteria the weak
equivalences and cofibrations should satisfy if (3) is desired.  As to (4),
in the present context it rests on a theorem of Jeffrey Smith that promises
to be valuable in the study of set-theoretically well-behaved Quillen model
categories; see Thm.~\ref{jeff}.  The main result of this paper,
\ref{auto}, is a meta-theorem to the following effect: if the ingredients
of a homotopy model category are given by suitable data so that they can be
interpreted in any topos, and over $\Set$ it does satisfy Quillen's axioms,
then it does so (functorially) in any topos.  There are mild conditions on
the syntax, and a single (annoying) set-theoretical one: cofibrations have
to be generated by a set.

We give six instances of the main result, some known, some new.  Perhaps it
is worthwhile to point out that the known cases had been obtained through
independent and fairly laborious methods.  The goal here is to invest
enough labor so that what \emph{is} tautological becomes visibly so.

Two corollaries may deserve attention: the existence of Quillen model
categories on simplicial objects in a topos with the class of cofibrations
smaller than all monos (Example~\ref{bkmono}), and existence of the
unbounded derived category of any Grothendieck abelian category
(Prop.~\ref{unbounded}).  This latter fact is certainly folklore, but the
proof via presentable model categories makes localization arguments easier
as well.

\bd{Motivation.} Perhaps it is useful to devote a few paragraphs to
sketching how a topos may be thought of as a geometric object, and why (and
when) it is worthwhile to do so.  This is independent of the technical
content matter of the paper, and readers acquainted with Grothendieckian
sheaf theory should skip to the next section.

Many types of geometric objects --- manifolds, orbifolds, algebraic spaces,
schemes --- are locally isomorphic to a fixed collection of distinguished
models, and the chief strength of sheaf theory is to make this precise and
to classify the extent to which such geometric objects may be
non-isomorphic globally.  Grothendieck observed that in many cases, the
natural notion of morphism between geometric objects leads to one and the
selfsame kind of adjunction between their associated categories of sheaves.
He suggested that topological (homotopical) invariants be thought of as
belonging to the abstract category of sheaves (the topos), functorially in
topos morphisms, rather than to the chosen representation in terms of a
site associated to the geometric object.  There are fruitful aspects of
this point of view.  An invariant constructed from a topos will have
avatars in many contexts, and several global properties of the category of
topoi --- for example, analogues of compactness, mapping spaces and
classifying spaces --- indeed make them similar to ``spaces''.  There are
also drawbacks, notably the lack (to date) of any ``cellular'' or
``skeletal'' or ``dimension-theoretic'' approach to an abstract topos,
though these are often the tools used to break up and understand a
geometric object.  (Instead, one may try to replace an unwieldy topos by a
more combinatorial structure with isomorphic invariants; for example, a
pro-simplicial set or sheaves on a simple Grothendieck topology, such as a
poset.)  But there is another, distinctly post-Grothendieckian approach,
which sees any category of set-valued sheaves as very similar in
\emph{some} formal properties to the category of sets.  Consequently the
homological (homotopical) algebra of sheaves of structures --- through
which one hopes to capture invariants of the topos --- ought to be similar
in \emph{some} formal aspects to plain homological (homotopical) algebra
with no Grothendieck topology in sight.  The onus is on making precise the
qualifier ``some''.  Let us believe that Quillen's axioms provide an
adequate calculus for homotopical manipulations.  The question I sought to
answer was: what formal properties should a homotopy theory satisfy if one
wishes to ``lift'' it into any category of sheaves as functorially as the
homological algebra of $R$-modules can be lifted?  Note that the topos
realm has proven to be a fertile ground for finding cohomological functors
(which ought to be thought of as just a part of homotopy theory) with
prescribed properties.  This note, however, is only concerned with a
foundational outlook.

\emph{Machinery.} The seeds of the techniques that provide an answer to the
motivating question were sown (with the exception of coherent logic) by
Grothendieck's school, but matured only later.  The
Makkai--Par\'e~\cite{makkpare} theory of accessibility is both a
simplification and a far-reaching extension of the notion of
$\pi$-accessibility (possibly the most technical part of SGA4).  See the
introduction to \cite{makkpare} for a full discussion of the prehistory.
The concept of locally presentable category was isolated by Gabriel and
Ulmer by pondering localizations of abelian categories and sheaf
reflections; cf.\ esp.\ Ulmer~\cite{ulmer}.  Both of these are, ultimately,
set-theoretical notions that allow one (in the present context) to treat
with ease arbitrary sites, even those lacking finiteness properties, points
etc.  As to the deployment of mathematical logic, exotic as it may seem in
a paper devoted to abstract homotopy theory, it is useful for two reasons.
First, it allows one to parametrize formal homotopy theories that have the
right functorial behavior without lapsing into vagueness as to what
constitutes a ``structure'' or a ``correct definition''.  (For example,
between the two ways of saying that a functor induces an equivalence of
categories --- namely, that it has a quasi-inverse, or that it is full,
faithful and essentially surjective, cf.~Ex.~\ref{sstack} --- the
difference is precisely that the second one proceeds within the language of
limits and colimits of graphs while the first one doesn't: it mentions
\emph{there exists a functor} that is a quasi-inverse$\dots$)  Secondly,
one has powerful theorems to the following effect: \emph{if} a mathematical
statement, formulated within the bounds of a specific logical syntax, holds
for structured sets \emph{then} it holds for like structured objects of any
topos.  These theorems are essential to the main Thm.~\ref{auto}.  Further
comparisons with existing work are given in the text in the form of remarks
which, I hope, do not distract the reader.

\emph{Omissions.} Theorem~\ref{auto} engulfs motivating examples \ref{joy}
and \ref{sstack}, but not \ref{siring}.  The reason is that the case of the
forgetful functor from simplicial algebraic theories to the underlying
objects, and many other such right adjoint situations, demand a slightly
more involved treatment, purely by virtue of the more complicated
description of cofibrations.  They are grouped together in the main theorem
of \cite{right}, part of which is concerned with setting up a logical
background that permits arguments similar to e.g.\ Cor.~\ref{condi2}.

Also, this paper is concerned solely with structures (in the sense of
footnote~\ref{sgafoot}) whose coefficients are constant, i.e.\ come from
the base topos $\Set$.  Thus, the results would apply to the category of
$R$-modules in a topos if $R$ is constant, but not to an arbitrary ringed
topos.  Many of the most interesting topos-theoretic Quillen model
categories take non-constant parameters; for example, Voevodsky's $I$-local
homotopy theory of simplicial objects, where $I$ is an ``interval object''
in the topos, or $G$-equivariant simplicial objects, where $G$ is a (not
necessarily constant) simplicial group.  The mixture of set-theoretic and
logical tools used in this paper is robust enough to apply to such
situations, but pursuing just those two threads is probably more important
than a bid at all-encompassing generality.

\bd{Acknowledgments.} I'd like to thank Michael Makkai and Ji\v{r}\'{\i}
Rosick\'y for making themselves accessible, Phil Hirschhorn and Charles
Rezk for being model correspondents, and Ieke Moerdijk and Jaap van Oosten
for being very logical.  My debt to Jeff Smith and Dan Kan is so great that
it needs a sentence of its own.

At the referee's suggestion, background material on formal homotopy theory
and on categorical logic is given in the body of the text, rather than in
an appendix.  This makes for quite a slow start in sections 1 and 2, but my
guess is that each reader will be able to skip a good half of the
preliminaries --- though possibly different halves.

\section{Presentable homotopy model categories}
\subsection{Homotopy model categories} In this article, the terms
\emph{Quillen model category} and \emph{homotopy model category} will both
refer to what Quillen~\cite{quil67} calls \emph{closed model category}.
(The reason for dropping the qualifier ``closed'' is that the interference
with \emph{closed category} is unfortunate, and the weaker ``non-closed''
axioms~\cite{quil67} will not be used.)  Good introductions to the subject
include Dwyer--Spali\'nski~\cite{model_cat}, Hovey~\cite{hovey} and
Goerss--Jardine~\cite{goerssjard2}.  Quillen's axioms~\cite{quil69} are:
\begin{itemize}
  \item[\bd{M1}:] $\C$ has finite limits and colimits. 
  \item[\bd{M2}:] If $f$ and $g$ are composable morphisms in $\C$, and if
    two of $f$, $g$ and $fg$ are weak equivalences, then so is the third. 
  \item[\bd{M3}:] A retract (in the category of morphisms of $\C$) of a
    fibration, cofibration or weak equivalence is respectively a
    fibration, cofibration or weak equivalence. 
  \item[\bd{M4}:] Given the commuting solid arrow diagram
    \[ \xymatrix{
        A \ar@{ >->}[d]_i \ar[r] & X \ar@{{}->>}[d]^p \\
        B \ar@{{}-->}[ur]^l \ar[r] & Y } \]
    with $i$ a cofibration and $p$ a fibration, if (i) $p$ or (ii) $i$ is
    a weak equivalence then a lifting $l$ exists making
    both triangles commute.  (One also says, ``$i$ has the left lifting
    property with respect to $p$'' or ``$p$ has the right lifting property
    with respect to $i$'' when an $l$ exists in every commutative square
    of this type.)
  \item[\bd{M5}:] Every morphism can be factored as (i) an acyclic
    cofibration followed by a fibration, and also as (ii) a cofibration
    followed by an acyclic fibration.
\end{itemize}
(Co)fibrations that are also weak equivalences will be called
\emph{acyclic} (rather than \emph{trivial} or \emph{aspherical}).

\subsection{Quillen adjoint pairs} Let $\M$, $\cN$ be homotopy model
categories and $\M\three{L}{\leri}{R}\cN$ an adjunction.  The following
three conditions are equivalent: (i) $L$ preserves cofibrations and $R$
preserves fibrations (ii) $L$ preserves cofibrations and acyclic
cofibrations (iii) $R$ preserves fibrations and acyclic fibrations.
$\HOMODEL$ is formed by the ``class'' of homotopy model categories with
these adjunctions as morphisms.\footnote{Since a typical Quillen model
category is not accessible, an extra Grothendieck universe is needed to
form $\HOMODEL$.  If one restricts to presentable homotopy model categories
as below, then $\HOMODEL$ lives in the same universe as $\TOPOI$; both are
subcategories of $\ACC$, the Makkai--Par\'e~\cite{makkpare} category of
accessible categories and functors.} By convention, the right adjoint gives
the direction of the arrow.  (2-categorical aspects are ignored; note, in
particular, that equivalence of homotopy model categories in Quillen's
sense does \emph{not} coincide with equivalence in $\HOMODEL$.)

\subsection{Transfinite composition} Let $\C$ be a category and $I$ a class
of morphisms of $\C$.  A morphism $F(0)\ra\colim F$ is said to be a 
\emph{transfinite composition} of arrows in $I$ if
\begin{itemize}
\item[$\bullet$] it is the part of a colimit cocone on $F:\alpha\ra\C$ from
$F(0)$ to the colimit.  Here $\alpha$ is an ordinal (thought of as an
ordered set, hence diagram) and $0\in\alpha$ its smallest element.
\item[$\bullet$] $F$ takes all successor arrows $\beta\prec\beta^+$ in
$\alpha$ into a morphism of $I$,
\item[$\bullet$] and $F$ is \emph{continuous}: for every limit ordinal 
$\beta\prec\alpha$, $F$ restricted to the diagram $\{\gamma\preceq\beta\}$
is a colimiting cocone in $\C$ on $F$ restricted to $\{\gamma\prec\beta\}$.
\end{itemize}

\begin{defn}
Let $\C$ be a cocomplete category, $I$ any class of morphisms of $\C$.
   \begin{itemize}
\item[$\bullet$] Close the class of all pushouts of $I$ under transfinite
composition.  This defines the class $\cell(I)$ of
\emph{relative $I$-cellular maps}.
\item[$\bullet$] The class $\cof(I)$ of $I$-cofibrations is defined as
follows: $X\llra{c}Y\in\cof(I)$ iff $c$ is a retract of an
$X\llra{r}Z\in\cell(I)$ in the category $X/\C$ of objects under $X$.
\item[$\bullet$] $I$-fibrations, or $I$-injectives, denoted $\inj(I)$, are
the morphisms with the right lifting property w.r.t.\ $I$; that is, such
that in any commutative square
\[ \xymatrix{
     \bullet \ar[d]_i \ar[r] & \bullet \ar[d]^p \\
     \bullet \ar@{{}..>}[ur] \ar[r] & \bullet } \]
with $i\in I$, $p\in\inj(I)$, a dotted lift making both triangles commute
exists.
   \end{itemize}
\end{defn}

\begin{remark}
Morally, a Quillen model category $\C$ is \emph{cofibrantly generated} if
(it is cocomplete and) the class of its cofibrations is of the form
$\cof(I)$ for some \emph{set} $I$, and its acyclic cofibrations are
$\cof(J)$ for some set $J$.  The fibrations (resp.\ acyclic fibrations)
must then be $\inj(J)$ ($\inj(I)$, resp).  The converse implication holds
if the small object argument applies; and that is the definition of
cofibrant generation chosen by Dwyer--Hirschhorn--Kan~\cite{dhk},
Hovey~\cite{hovey}.  If a Quillen model category is locally presentable, as
will always be in this article, it is cocomplete and the small object
argument (and more) applies to any set of maps; hence, for the purposes of
this paper, cofibrant generation is the ``moral'' referred to above.
\end{remark}

\subsection{The transfinite small object argument} The inventor's account
(with a mild inaccuracy related to not requiring a certain cardinal to be
regular) is Bousfield~\cite{bous76}.  The proof will be quite abbreviated
since so many versions are already in the literature; see e.g.\
Hovey~\cite{hovey} (which builds on Hirschhorn~\cite{hirsch} and
Dwyer--Hirschhorn--Kan~\cite{dhk}) for full details.  For an introduction
to the lore of locally presentable and accessible categories, see \AdRos\
or Borceux~\cite{borceux} vol.II.

\begin{prop} \label{trasmall}
Let $\C$ be a locally presentable category and $I$ a \emph{set} of 
morphisms of $\C$.
     \begin{itemize}
\item[--] Every $m\in\mor\C$ can be factored as $fc$ with $c\in\cell(I)$,
$f\in\inj(I)$.  The factorization is not unique, but can be performed
functorially.
\item[--] $\cof(I)$ is exactly the class of morphisms having the left
lifting property w.r.t.\ all elements of $\inj(I)$, and $\inj(I)$ is exactly
the class of morphisms with the right lifting property w.r.t.\ all elements
of $\cof(I)$.
     \end{itemize}
\end{prop}

\begin{proof}
Define a functor $\Mor(\C)\llra{F}\C$ and natural transformations
$\dom\llra{\Theta}F$, $F\llra{\Psi}\codom$ (these being the domain resp.\
codomain functors) factoring the identity $\Mor(\C)\ra\Mor(\C)$ as follows.
For a map $X\llra{m}Y$ of $\C$, let $S$ be the set of all commutative 
diagrams
\[ \xymatrix{    A \ar[d]_i \ar[r] & X \ar[d]^m \\
                 B \ar[r]          & Y } \]
such that $i\in I$, and $A_S\ra B_S$ the coproduct of all these arrows $i$;
$A_S$ comes with a natural map to $X$.  The value of $F$ at $m$ is the
pushout of
\[ \xymatrix{ A_S \ar[d] \ar[r] & X  \\
              B_S } \]
with its natural induced map $\Theta(m)$ (resp.\ $\Psi(m)$) from $X$
(resp.\ to $Y$).  Beginning with $F_0:=F$, $\Theta_0:=\Theta$,
$\Psi_0:=\Psi$, build by transfinite induction a diagram of functors
$\Mor(\C)\llra{F_\alpha}\C$ and compatible natural transformations
$\dom\llra{\Theta_\alpha}F_\alpha$, $F_\alpha\llra{\Psi_\alpha}\codom$.
For a successor ordinal $\alpha+1$, $F_{\alpha+1}(-):=F(\Psi_\alpha(-))$, 
$\Psi_{\alpha+1}(-):=\Psi(\Psi_\alpha(-))$,
$\Theta_{\alpha+1}:=\Theta(\Psi_\alpha(-))\circ\Theta_\alpha$.  For a limit
ordinal $\alpha$, $F_\alpha$ resp.\ $\Psi_\alpha$, $\Theta_\alpha$ are
colimits along the chain $\beta\prec\alpha$.  Let $\kappa$ be a regular
cardinal greater than the rank of presentability of the domain of any arrow
in $I$.  The requisite factorization of $m$ is 
$X\llra{\Theta_\kappa(m)}F_\kappa(m)\llra{\Psi_\kappa(m)}Y$. 
($\Theta_\kappa(m)$ is a transfinite composition of coproducts of pushouts
of arrows from $I$, but a coproduct of arrows is a transfinite composition
of pushouts, and a transfinite composition of transfinite compositions is
one such again.)  The second part of the claim follows by factoring any
$m\in\cof(I)$ as $fc$ with $f\in\inj(I)$ and $c\in\cell(I)$, and noting
that $m$ has the left lifting property w.r.t.\ $\inj(I)$, in particular,
w.r.t.\ $f$, which entails that it is a retract of $c$ of the type claimed.
\end{proof}

\begin{remark}
For the proof to work, it is enough for $\C$ to be cocomplete and for the
domains $X$ of the maps in $I$ to be such that $\hom_\C(X,-)$ commutes with
transfinite compositions of morphisms from $\cell(I)$ provided they are
``long enough''.  The assumption that $\C$ is locally presentable is much
stronger; for example, it implies the existence of a set of dense
generators.  There exist cocomplete categories other than locally
presentable ones where every object $X$ has a \emph{rank}, i.e.\
$\hom(X,-)$ commutes with all $\kappa$-filtered colimits for some $\kappa$
depending on $X$ (take e.g.\ free cocompletions of certain \emph{large}
categories).  Finally, even such categories as topological spaces or
various topological spectra, where not every object has a rank, may possess
\emph{certain} sets of morphisms $I$ which ``permit the small object
argument'', in the terminology of Dwyer--Hirschhorn--Kan~\cite{dhk},
Hovey~\cite{hovey}.  In our context of algebraic models for homotopy types,
however, local presentability is a convenient ground assumption.
\end{remark}

\begin{defn}  (the solution set conditions) \\  \label{solset}
Let $\C$ be a category, $\tW$ a class of morphisms, $m\in\mor\C$.  Say
that $\tW$ \emph{satisfies the solution set condition at} $m$ if there
exists a subset $\tW_m$ of $\tW$ such that any
commutative square
\[\xymatrix{ \bullet \ar[r] \ar[d]_m & \bullet \ar[d]^w \\
                    \bullet \ar[r]   & \bullet                } \]
with $w\in\tW$ allows to be factorized by a commutative diagram
\[\xymatrix{ \bullet \ar[r] \ar[d]_i & \star \ar[d]^{w_m} \ar[r]
                                                 & \bullet \ar[d]^w \\
             \bullet \ar[r] & \star \ar[r] & \bullet       } \]
with $w_m\in\tW_m$.  Let $I$ be a class of morphisms; $\tW$ satisfies the
solution set condition at $I$ if it satisfies it at each $m\in I$.  If the
solution set condition is satisfied for every $m\in\mor\C$, say simply
that $\tW$ satisfies the solution set condition.
\end{defn}

\begin{remark}  \label{freyd}
Let $\Mor(\C)$ be the category of morphisms of $\C$ (maps in $\Mor(\C)$
are commutative squares in $\C$), and $\Mor(\tW)$ the full subcategory of 
$\Mor(\C)$ whose objects belong to $\tW$.  Def.~\ref{solset} is Freyd's
solution set condition for the inclusion functor $\Mor(\tW)\inc\Mor(\C)$.
(Even if $\tW$ is a subcategory of $\C$, the inclusion $\tW\inc\C$ will
not be considered in this context, so the terminology of Def.~\ref{solset}
should result in no confusion.)
\end{remark}

The following theorem was announced by Jeffrey Smith at the 1998 Barcelona
conference in Algebraic Topology.  It greatly amplifies and simplifies
results of Goerss and Jardine~\cite{goerssjard} and the
author~\cite{tezis}.  I am indebted to him for explaining his proof, and
for his permission to reproduce it here. 

\begin{thm} \label{jeff} \textup{(J.~Smith)} \\
Let $\C$ be a locally presentable category, $\tW$ a subcategory, and $I$ a
set of morphisms of $\C$.  Suppose they satisfy the criteria:
\begin{itemize}
\item[\bd{c0}] $\tW$ is closed under retracts and has the 2-of-3 property
(Quillen's axiom \bd{M2}).
\item[\bd{c1}] $\inj(I)\subseteq\tW$.
\item[\bd{c2}] The class $\cof(I)\cap\tW$ is closed under transfinite
composition and under pushout.
\item[\bd{c3}] $\tW$ satisfies the solution set condition at $I$.
\end{itemize}
Then setting weak equivalences:=$\tW$, cofibrations:=$\cof(I)$ and
fibrations:=$\inj(\cof(I)\cap\tW)$, one obtains a cofibrantly generated
Quillen model structure on $\C$.
\end{thm}

\begin{proof} (J.~Smith)
The strategy is to exhibit a set $J$ of morphisms such that
$\cof(J)=\cof(I)\cap\tW$.  From there, Quillen's axioms follow in a
well-known way.  Two small object arguments yield the factorization axiom
\bd{M5}; the part of \bd{M4} that is not the definition follows from 
\bd{M5}, \bd{M2} and the retract argument; \bd{M3} holds by the definition
of (co)fibrations and \bd{c0}; finally, a locally presentable category is
complete and cocomplete, so \bd{M1} is satisfied. 

$J$ itself will be constructed in two steps.  Lemma~\ref{dense} shows that
if a collection $J$ of morphisms is ``dense'' between $I$ and $\tW$, then
$\cof(J)=\cof(I)\cap\tW$.  Lemma~\ref{small}, using \bd{c3}, constructs
such a $J$ that is only a set.

\begin{lemma} \label{dense}
Let $J\subseteq\cof(I)\cap\tW$ be a collection (set or possibly proper
class) of maps in $\C$ such that for any commutative square
\[ \xymatrix{ \bullet \ar[d]_i \ar[r] & \bullet \ar[d]^w \\
              \bullet \ar[r]          & \bullet                } \]
with $i\in I$, $w\in\tW$ there exists $j\in J$ that factors it:
\[\xymatrix{ \bullet \ar[r] \ar[d]_i & \star \ar[d]^j \ar[r]
                                                 & \bullet \ar[d]^w \\
             \bullet \ar[r] & \star \ar[r] & \bullet           } \]
Then any $f\in\tW$ can be factored as $hg$ with $g\in\cell(J)$,
$h\in\inj(I)$.
\end{lemma}

\bd{Corollary.}
Under the assumptions of the previous lemma, $\cof(J)=\cof(I)\cap\tW$.

\emph{Proof of the corollary.}
$\cof(J)$ is the saturation of $J$ under pushout, transfinite composition
and retracts, $J\subseteq\cof(I)\cap\tW$ and $\cof(I)\cap\tW$ is supposed
to be closed under these operations, so $\cof(J)\subseteq\cof(I)\cap\tW$.
Conversely, consider any $f\in\cof(I)\cap\tW$ and write it $f=hg$ as above.
Since $f\in\cof(I)$ and $h\in\inj(I)$, $f$ is a retract of $g$ (in the
category of objects under the domain of $f$).  So $f\in\cof(J)$.

\emph{Proof of Lemma~\ref{dense}.}
This is rather like the ordinary small object argument, save that one glues
on the ``interpolating'' maps $J$ instead of the $I$.  More precisely, we
wish to build by transfinite induction on $\lambda$ certain factorizations
\[ X=:P_0\ra P_1\ra\dots\ra P_\alpha\ra P_{\alpha+1}\ra\dots\ra P_\lambda
\llra{h_\lambda} Y \]
of $f$ such that (this is the induction hypothesis) the diagram 
$P_0\ra\dots\ra P_\lambda$ is a continuous composition of maps belonging to
$\cell(J)$.  Thence the composite itself will belong to $\cell(J)$.  Since
$\cell(J)\subseteq\tW$ and $f\in\tW$, the 2-of-3 property of $\tW$ implies
that $h_\lambda\in\tW$. \\
Set $P_0:=X$, $h_0:=f$.  At a successor stage, let $S_\lambda$ be the set
of all commutative squares
\[\xymatrix{ A \ar[r] \ar[d]_i & P_\lambda \ar[d]^{h_\lambda} \\
                    B \ar[r]   & Y                          } \]
with $i\in I$.  The density assumption on $J$ means the existence of a
factorization
\[\xymatrix{ \bullet \ar[r] \ar[d]_i & A_s \ar[d]^{j_s} \ar[r]^{t_s}
                                         & P_\lambda \ar[d]^{h_\lambda} \\
             \bullet \ar[r] & B_s \ar[r] & Y        } \]
with $j_s\in J$, for each square $s\in S_\lambda$.  Let $P_{\lambda+1}$ be
the pushout
\[\xymatrix{\coprod A_s \ar[r]\ar[d]_{\coprod j_s} &  P_\lambda \ar[d] \\
            \coprod B_s \ar[r]   & P_{\lambda+1}  } \]
along the canonical
$\coprod_{s\in S_\lambda}A_s\xra{\{t_s\;|\;s\in S_\lambda\}}P_\lambda$.
Let $h_{\lambda+1}$ be the canonical pushout corner map from 
$P_{\lambda+1}$ to $Y$.  The connecting map $P_\lambda\ra P_{\lambda+1}$ is
a pushout of coproducts of morphisms from $J$.  But any coproduct of maps
is a transfinite composition (starting from the coproduct of the domains),
so the connecting map belongs to $\cell(J)$. \\
At a limit ordinal $\lambda$, $X\ra P_\lambda$ is the colimit of the
diagram $X\ra\dots\ra P_\alpha$, $\alpha\prec\lambda$.

Let now $\kappa$ be a regular cardinal exceeding the rank of presentability
of all the objects that occur as domains of maps in $I$.  The required
factorization of $f$ is $X\llra{g}P_\kappa\llra{h}Y$.

Indeed, consider any lifting problem
\[\xymatrix{ A \ar[r]^a \ar[d]_i & P_\kappa \ar[d]^h  \\
             B \ar[r]   & Y                          } \]
with $i\in I$.  Since $\kappa$ is regular, the diagram $X\ra\dots\ra
P_\kappa$ is $\kappa$-filtered, and since $\hom(A,-)$ commutes with
$\kappa$-filtered colimits by assumption, $a$ factors through a prior stage
$A\ra P_\lambda\ra P_\kappa$.  If the lifting problem
\[\xymatrix{ A \ar[r] \ar[dd]_i & P_\lambda \ar[d]  \\
                                & P_\kappa  \ar[d]^h  \\
             B \ar[r]   & Y                          } \]
is indexed by $s\in S_\lambda$, the solution to the original one is the
bottom composite
\[\xymatrix{
A \ar[r]\ar[d]_i & \star \ar[d]^{j_s}\ar[r]
            & \bigstar \ar[d]^{\coprod j_s} \ar[r] & P_\lambda \ar[d] \\
B \ar[r]  & \star \ar[r] & \bigstar \ar[r] & P_{\lambda+1} \ar[r]
& P_\kappa.  } \]

\begin{lemma} \label{small}
There exists a set $J$ with the property required in Lemma~\ref{dense}.
\end{lemma}

Indeed, consider the set of all morphisms (in the category of arrows) from
$i\in I$ to the solution set $\tW_i$:
\[\xymatrix{ \bullet \ar[r] \ar[d]_i & X \ar[d]^{w_i\in\tW_i}  \\
             \bullet \ar[r]          & Y                        } \]
form the pushout $P$ and the canonical corner map $c$
\[\xymatrix{ \bullet \ar[r]\ar[d]_i & X \ar[rdd]^{w_i} \ar[d]_{i'} \\
               \bullet \ar[r] \ar[drr] & P \ar[rd]|c \\
                                                && Y            } \]
and factor $c$ as $P\llra{p}Q\llra{q}Y$ with $p\in\cell(I)$, $q\in\inj(I)$.
Set $j:=pi'$.  $J$ is the set of such $j$ (one for each morphism from
$i\in I$ to $\tW_i$).  Indeed, $i'\in\cell(I)$ and $p\in\cell(I)$, so
$j\in\cof(I)$.  $q\in\inj(I)\subseteq\tW$ by \bd{c1}, so $w_0=qpi'\in\tW$
and 2-of-3 imply $j\in\tW$.  Finally, any morphism from $I$ to $\tW$
\[ \xymatrix{ \bullet \ar[d]_i \ar[r] & \bullet \ar[d]^w \\
              \bullet \ar[r]          & \bullet                } \]
allows to be factored, using \bd{c3}, as
\[ \xymatrix{ \bullet \ar[d]_i \ar[r] &  X \ar@{=}[r] \ar[d]_j
                             & X \ar[d]^{w_i} \ar[r] & \bullet \ar[d]^w \\
\bullet \ar[r] & Q \ar[r]^q & Y \ar[r] & \bullet. } \]
This completes the proof of lemma~\ref{small} and of Jeff Smith's theorem.
\end{proof}

\begin{remark}
There are numerous variations on the above themes, notably
Dwyer-Hirschhorn-Kan~\cite{dhk}, Hirschhorn~\cite{hirsch} and
Stanley~\cite{stanley}, but J.~Smith's theorem seems to be the first to
identify Freyd's solution set condition as the ``culprit'' for there being
a set of generating acyclic cofibrations, provided a set of generating
cofibrations is known to exist.
\end{remark}

\begin{remark}
By virtue of the way the set $J$ is found in lemma~\ref{small},
lemma~\ref{dense} would go through by assuming merely that the domains of
$I$ are small w.r.t.\ long enough transfinite $\cell(I)$-compositions.  To
prove \bd{M5} however, one also needs to do a small object argument on $J$.
Since the solution set is non-canonical, it is best to assume that for
every object $X$, $\hom(X,-)$ commutes with $\kappa$-filtered sequential
colimits for some $\kappa$; and this is guaranteed by the assumption that
$\C$ is locally presentable.
\end{remark}

There is a notable type of Quillen model categories whose cofibrations are
identifiable explicitly as the monomorphisms.  In all such cases I am aware
of, the following proposition applies to show that they are generated by a
set.  The proof really goes back to Grothendieck; see also
Barr~\cite{barr88}.

\begin{prop}   \label{mono}
Let $\C$ be a category; write $\mono$ for its class of monomorphisms.
Suppose
  \begin{itemize}
\item[(i)] $\C$ is locally presentable.
\item[(ii)] Subobjects have effective unions in $\C$.  That is,
\[ \xymatrixcolsep{1pc} \xymatrixrowsep{1pc} \xymatrix{
  A \cap B \ar@{ >->}[rrr] \ar@{ >->}[ddd]
                            &&& A \ar@{ >->}[ddd] \ar[ddl]_a \\ \\
          && A \cup B \ar@{ >->}[dr]^m \\
   B \ar@{ >->}[rrr] \ar[rru]^b &&& X } \]
given any two subobjects $A,B$ of an object $X$, form their intersection
$A\cap B=A\times_X B$ and their pushout $A\cup B$ over their intersection;
the induced maps $a,b,m$ are to be monomorphisms (whence $A\cup B$ really
is the supremum of $A$ and $B$ in the subobject lattice of $X$, i.e.\ their
union).
\item[(iii)] $\mono$ is closed under transfinite composition.
  \end{itemize}
Then $\mono=\cell(I)=\cof(I)$ for some set $I\subset\mono$.
\end{prop}

\begin{remark}
Any $AB5$ abelian category and any elementary topos satisfies (ii)-(iii). 
(i) restricts them to Grothendieck abelian categories (see 
Prop.~\ref{grolo}) resp.\ Grothendieck topoi.  If a category satisfies 
(i)-(ii)-(iii), so do all its over-, under-, and diagram categories, and
left exact localizations.
\end{remark}

\begin{proof}
Let $\G$ be a set of strong generators for $\C$, so that the functors
$\hom(G,-)$, $G\in\G$, collectively reflect isomorphisms.  (Any locally
presentable category has such $\G$.)  Let $\Q$ be the set of (isomorphism
types of) regular quotients of these generators; finally, let $I$ be the
set of all (isomorphism types of) subobjects of members of $\Q$.  Then
$\mono=\cell(I)$ (a fortiori $\mono=\cof(I)$, since monos are closed under
retract).

Argue by contradiction.  Suppose $X\llra{m}Y$ is a mono but
$m\not\in\cell(I)$.  By transfinite induction, we will build a chain
$X:=P_0\ra P_1\ra P_2\ra\dots\ra P_\lambda\ra\dots\ra Y$ of subobjects of
$Y$ that is (i) properly increasing (ii) satisfies $X\ra
P_\lambda\in\cell(I)$ for $0\prec\lambda$.  This contradicts that $\C$ is
well-powered.

Set $P_0:=X\llra{m}Y$.  At a successor stage, find $G\llra{g}Y$, $G\in\G$,
that does not factor through $P_\lambda\ra Y$.  Factor $g$ as a regular epi
followed by a mono: $G\epi Q\mon Y$.  Form an effective subobject union
diagram as above:
\[ \xymatrixcolsep{1pc} \xymatrixrowsep{1pc} \xymatrix{
 & Q\cap P_\lambda \ar@{ >->}[rrr] \ar@{ >->}[ddd]
                        &&& P_\lambda \ar@{ >->}[ddd] \ar[ddl]^a \\ \\
          &&& Q\cup P_\lambda \ar@{ >->}[dr] \\
 G \ar@{->>}[r] &  Q \ar@{ >->}[rrr] \ar[rru] &&& Y } \]
and define $P_{\lambda+1}:=Q\cup P_\lambda$.  Note $a\in\cell(I)$.
$P_{\lambda+1}$ is bigger than $P_\lambda$ since $g$ factors through it, so
the induction hypotheses are satisfied.

At a limit ordinal $\lambda$, set
$P_\lambda:=\colim_{\alpha\prec\lambda}P_\alpha$ and use assumption (iii).
\end{proof}

\subsection{Solution sets vs.\ accessibility}
We recall results on the second named set-theoretic constraint.  In the
context of classes of morphisms, it is stronger than \ref{solset} but much
better behaved under (2-)categorical operations.  All in all, it may be
easier to check for accessibility than to solve for a solution set.  The
interaction with homotopy theory (or, rather, Quillen model categories
whose underlying category is locally presentable) is solely through the
solution set condition for weak equivalences, and their closure under
retracts (see Prop.~\ref{retrac}).  The reader may wish to skip to the next
section and refer back only as necessary.

\begin{defn} (accessibility of a class of maps) \\
Let $\C$ be a locally presentable category, $\tW$ a class of morphisms of
$\C$.  Let $\Mor(\tW)$ and $\Mor(\C)$ be as in remark~\ref{freyd}.  Say
that \emph{$\tW$ is an accessible class of maps}\footnote{The last part of
remark~\ref{freyd} applies here as well.} if $\Mor(\tW)$ is an accessible
category, and it is closed in $\Mor(\C)$ under $\kappa$-filtered colimits
for some $\kappa$.
\end{defn}

In the terminology of \AdRos, $\Mor(\tW)$ is an accessibly embedded,
accessible subcategory of $\Mor(\C)$.  We refer to \AdRos,
Borceux~\cite{borceux} vol.II or Makkai--Par\'e~\cite{makkpare} for further
background.  Note that accessibility of a class of maps lacks meaning
unless the ground category is locally presentable (or accessible, at
least).  In what follows, fix a locally presentable $\C$ and class of maps
$\tW$.

\begin{prop}  \label{accsolset}
If $\tW$ is accessible, it satisfies the solution set condition.
\end{prop}

\begin{proof}
More generally, any accessible functor satisfies the solution set condition
(at every object); see \AdRos~Corollary 2.45.
\end{proof}

\begin{remark}
The distinction between being an accessible class and satisfying the
solution set condition is subtle.  (These notions have the obvious meaning
for any class of objects in a locally presentable category $\K$, and the
statements about to be quoted apply to this more general case.)  A theorem
due to H.~Hu and M.~Makkai~\cite{humakk} asserts that a class of objects
closed in $\K$ under $\kappa$-filtered colimits (for some $\kappa$) is 
accessible iff it satisfies the solution set condition.  J.~Rosick\'y and
W.~Tholen~\cite{rostho} prove that the set-theoretical statement known as
Vop\v{e}nka's Principle (a so-called \emph{large cardinal axiom}) implies
that any class of objects satisfying the solution set condition will be
accessible.  This diminishes the chances of finding, without the aid of
axioms external to ZFC, any class of objects that is not accessible but
satisfies the solution set condition.  Such a counterexample would disprove
Vop\v{e}nka's Principle from ZFC, and current set theoretical intuition is
that this is unlikely.
\end{remark}

\begin{remark}
One reason to cherish $\ACC$, the (very large) category of accessible
categories and functors is the Limit theorem of Makkai and Par\'e
\cite{makkpare}: $\ACC$ is closed in the (very large) category of (large) 
categories under weighted pseudo-limits.  As a corollary, $\ACC$ is stable
under forming under-, over-, comma-, diagram categories$\dots$; the 
intersection of a set of accessible subcategories is accessible$\dots$ and
more.  A category is locally presentable iff it is accessible and
cocomplete; as further corollaries, the category of algebras (over an
accessible monad), the category of coalgebras (over an accessible comonad),
the category of [commutative] monoids (w.r.t.\ an accessible [symmetric]
monoidal structure) is locally presentable if the ground category is so.
Here we single out one corollary that ties in quite directly with the 
subject matter.
\end{remark}

\begin{prop}   \label{inv}
Let $\Mor(\A)\llra{F}\Mor(\B)$ be an accessible functor, and $\tW$ an
accessible class in $\B$.  $F^{-1}(\tW)$ (i.e.\ the class of morphisms in
$\A$ taken into $\tW$ by $F$) is an accessible class.
\end{prop}

\begin{proof}
The full inverse image of a full, accessible subcategory by an accessible
functor is again accessible; see \AdRos~Remark~2.50.
\end{proof}

In contrast to the previous ones, the next fact is elementary. 

\begin{prop}  \label{retrac}
Any accessible class $\tW$ of maps is closed under retracts (in the
category of morphisms).
\end{prop}

\begin{proof}
In fact, any full subcategory of cocomplete category closed under 
$\kappa$-filtered colimits for some $\kappa$ must be closed under retracts. 
This is since a retract is a colimit on an ``idempotent loop'' diagram
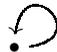
\[ \def\objectstyle{\scriptstyle}\xymatrix{\bullet\ar@(r,u)[]} \]
which is $\infty$-filtered, ie.\ $\kappa$-filtered for every $\kappa$.
\end{proof}

Here is a way to make homotopical machinery accessible.  By virtue of Kan's
combinatorial description of weak equivalences, the class of weak
equivalences in $\SSet$ is accessible.\footnote{See Example~\ref{sset}. The
claim does not appear to be a formal consequence of the definition /
theorem that a map in $\SSet$ is a weak equivalence iff its geometric
realization is.  There is no convenient category of topological spaces
known that is also locally presentable.} Let now $\C$ be a locally
presentable category, to be made into a Quillen model category.  Suppose
there exists a detection functor: $\C\llra{F}\SSet$ such that $w\in\mor\C$
is a weak equivalence iff $F(w)$ is one.  If $F$ is accessible, that is to
say, it preserves $\kappa$-filtered colimits for some $\kappa$, then weak
equivalences form an accessible class in $\C$, and \bd{c0} and \bd{c3} of
Theorem~\ref{jeff} are automatic.  J.~Smith conjectures that every
cofibrantly generated model category whose underlying category is locally
presentable, arises this way.

\section{Defining and ``sheafifying'' homotopy model theories} What is
common to simplicial objects, simplicial rings and categories is that they
are structures definable in terms of finite limits.  On the intuitive
level, this notion has been around since the beginning of category
theory\footnote{SGA4, Tome 1, Expose 1, 2.9 speaks thus: ``Soit $\gamma$
une esp\`ece de structure alg\'ebrique `d\'efinie par limites projectives
finies'. (Le lecteur est pri\'e de donner un sens math\'ematique \`a la
phrase pr\'ec\'edente.\label{sgafoot} Notons seulement que les structures
de groupes, groupes ab\'eliens, anneux, modules, etc$\dots$ sont de telles
structures.)''} and has been given several equivalent formal definitions
since then.

\subsection{Logics} A \emph{language} for (first-order, many-sorted) logic
consists of a set of \emph{types} (also called \emph{sorts}), function
symbols, relations symbols, logical connectives $\es,\vagy,\implies,\neg$
and quantifiers $\forall,\exists$.  For convenience, we also add the
symbols $\top$, $\bot$ and $\exists!$ for the logical constants ``true'',
``false'' and the quantifier ``there exists precisely one''.\footnote{One
could, in several ways, start with a smaller set of logical symbols as
basic and express the rest in terms of them, but there's no reason for such
parsimony for us.} There's to be a set of (dummy) \emph{variables}, each
assigned one of the types (and thought of as running over the set of things
of that type).  Similarly, each function symbol is supposed to come with
its \emph{arity} $\lc t_1,t_2,\dots,t_n\rc\ra t$, meaning it takes an
$n$-tuple of things, the $\text{i}^\th$ of which is type $t_i$, and turns
them into a thing of type $t$.  Analogously for relation symbols.  
\emph{Constants} may be thought of as 0-ary functions.  Refer to any
textbook for the precise formation rules of \emph{well-formed formulas}. A
formula containing no free variables (i.e.\ such that any variable falls
under the scope of a quantifier) is called a \emph{sentence}.  An axiom
system is a logical language together with a set of sentences (axioms).

\begin{example}
In the most natural way of formalizing the notion of \emph{vector space},
there exist two types (scalars and vectors), functions of ``scalar
multiplication'' (of arity $\lc$scalar,vector$\rc\ra$ vector), vector 
addition$\dots$, symbols for three constants: the zero vector and the 
scalars 0 and 1, numerous axioms, and no relation symbols other than 
equality.  (Since only things of the same type can be equal, and every 
relation symbol is to have its arity, there's supposed to be a separate $=$ 
for comparing two scalars, and one for comparing two vectors.)
\end{example}

An \emph{interpretation} of a language $\L$ (in the category $\Set$) is the
assignment of a set $T_i$ to each type $t_i$ of $\L$, an actual function
$T_1\times T_2\times\dots\times T_n\ra T$ for each function symbol of the
respective arity, and a subset $R$ of $T_1\times T_2\times\dots\times T_n$
for each relation symbol $r$ of arity $\lc t_1,t_2,\dots,t_n\rc$.  It is a
\emph{model} if the axioms are satisfied.

\subsubsection{The category of models} Models for a fixed axiom system
form a category, in fact a full subcategory of the category of 
interpretations of the language.  Let $M, N$ be models, with underlying
sets for types $T_i^M, T_i^N$, resp.  A morphism from $M$ to $N$ is given
by a function $T_i^M\llra{f_i}T_i^N$ for each type such that they commute
with the interpretations of the function symbols:
\[ \xymatrix{
T_1^M\times T_2^M\times\dots\times T_n^M \ar[d]^{\times\{f_i\}}
\ar[r] & T^M \ar[d]^{f} \\
T_1^N\times T_2^N\times\dots\times T_n^N \ar[r] & T^N
} \]
and for each relation symbol $r$, $R^M$ is to be taken into a subset of
$R^N$ by $T_1^M\times T_2^M\times\dots\times T_n^M\llra{\times\{f_i\}}
T_1^N\times T_2^N\times\dots\times T_n^N$.

\subsubsection{Interpreting logic in a topos} Let now $\E$ be a topos.  An
interpretation of a language in $\E$ is the assignment of an object $T_i$
to each type $t_i$ there is, a morphism $T_1\times T_2\times\dots\times
T_n\ra T$ for each function symbol, and a subobject $R$ of $T_1\times
T_2\times\dots\times T_n$ for each relation symbol of the respective arity.  
Using a calculus of subobjects --- intersections, unions,
pseudo-complements, direct and inverse images of subobjects --- and via
induction on the complexity of formulas, every sentence of $\L$ gets
assigned a truth value of ``true'' or ``false''.  The rules reduce to the
usual ones when $\E=\Set$; however, it is not the case that two formulas
that are logically equivalent in $\Set$ (i.e.\ always evaluate to the same
truth value) will always be logically equivalent in any topos.  The reader
is referred to the textbooks of MacLane--Moerdijk~\cite{m&m} or
Borceux~\cite{borceux}~vol.III for details.  The category of models in $\E$
of an axiom system can now be defined by analogy with the case of $\Set$.

In a \emph{fragment of logic} only expressions of some specific form are
permitted.  In this paper, the emphasis is on the following two fragments.

\subsection{Cartesian logic} A language for cartesian logic\footnote{See
Johnstone~\cite{john82} p.180 and p.221 for the origin of this term.} is a
(many-sorted, first-order) language containing \emph{no relation symbols
other than equality} and making no use of the logical symbols
$\vagy,\neg,\exists,\bot$.  Moreover, in a cartesian axiom system, only
sentences of the following kind are allowed:
\[     \forall\vec{x}\bigl(\Phi\implies\exists!\vec{y}(\Psi)\bigr)  \]
where $\vec{x}$, $\vec{y}$ are shorthand for a string of variables, and
$\Phi$ and $\Psi$ must be finite conjunctions of expressions of the form
$\bd{t}_1=\bd{t}_2$ where $\bd{t}_1$, $\bd{t}_2$ are \emph{terms}, that
is, meaningful combinations of constants, dummy variables and function
symbols. 

\begin{example}
Here's how to say \emph{groupoid} in this language.  There are two sorts of
things, objects $A_0$ and arrows $A_1$.  There are the usual functions of
``source'' \bd{s}, ``target'' \bd{t}, ``identity'' \bd{i}, and ``inverse'';
their commutation relations (to conform to the standard that implication
has to be used precisely once$\dots$) can be expressed in the form
\[      \forall a\bigl(\top\implies \bd{s}(\bd{i}(a))=a\bigr)      \]
and so on (here $\top$ is the True, and $a$ is a variable of the type 
\emph{object}).  To express composition and its associativity, one has to
introduce the auxiliary type $P$ for the ``set of'' composable pairs
of arrows.  There are projections $P\llra{\bd{p}_1}A_1$,
$P\llra{\bd{p}_2}A_1$; composition is a function $P\ra A_1$.  Another
typical axiom is
\[ \forall a_1,a_2\bigl(\bd{t}(a_1)=\bd{s}(a_2)\implies
\exists!p(\bd{p}_1(p)=a_1 \es \bd{p}_2(p)=a_2)  \bigr) \]
where $p$ is a variable of type $P$.
\end{example}

In fact, groupoids are (qualitatively) the most complicated type of example
that can occur.  Every axiom
$\forall\vec{x}\bigl(\Phi\implies\exists!\vec{y}(\Psi)\bigr)$ can be seen
as expressing the existence of a (vector-valued) function with domain
specified by $\Phi$ and codomain specified by $\Psi$, and asserting some
equational conditions about components of that vector.  So every cartesian
structure is an equational theory of partial algebras, and actually of a
special class of partial algebras: the domain of the partial operations has
to be specified by a conjunction of equational conditions between total
operations (morally, ``by pullbacks'').

Cartesian logic provides a solution to the exercise of SGA4 recalled in
footnote~\ref{sgafoot}: ``define \emph{definable in terms of finite
limits}''.  M.~Coste's \emph{lim-logic} (see Coste~\cite{coste}), Lawvere's
\emph{cartesian theories} (see Barwise~\cite{barwise} A.8,
Barr--Wells~\cite{barrwells}), Ehresmann's \emph{finite limit sketches}
(see \AdRos, Barr--Wells~\cite{barrwells}, Borceux~\cite{borceux} vol.III)
and P.~Freyd's \emph{essentially algebraic theories} (see \AdRos) have the
same expressive strength in a topos as cartesian logic does; in fact, any
two of them are functorially interdefinable.

The following proposition compresses all that will be needed of finite
limit$=$cartesian structures below.

\begin{prop}   \label{flimfact}
Given a structure $S$ defined in terms of finite limits, one can associate
to a topos $\E$ the category of $S$-structures (or ``models of the axiom
system $S$'') in $\E$, to be denoted $\Mod_S(\E)$.\footnote{Let $(\C,J)$
be a site; $\Mod_S(\sh(\C,J))$ is equivalent to the category of
$\Mod_S(\Set)$-valued sheaves on $(\C,J)$.  E.g.\ a sheaf of abelian groups
is an abelian group object in sheaves.  Finite limit definable structures
form the largest class of first-order structures to which this extends.}
\begin{itemize}
   \item[(i)] $\Mod_S(\E)$ is a locally presentable category.
   \item[(ii)] A topos morphism $\E\llra{f}\F$ (i.e.\ adjoint pair
$f^*\dashv f_*$ such that $f^*$ preserves finite limits) induces an
adjunction between $\Mod_S(\E)$ and $\Mod_S(\F)$.
   \item[(iii)] For any small category $\D$, $\Mod_S(\E)^\D$ is canonically
isomorphic to $\Mod_S(\E^\D)$.
   \item[(iv)] A morphism of $S$-structures is itself definable in
terms of finite limits; more precisely, there exists a finite limit
structure $\Mor(S)$ such that $\Mor(\Mod_S(\E))$ is canonically equivalent
to $\Mod_{\Mor(S)}(\E)$, for any $\E$.
\end{itemize}
\end{prop}

\begin{example}
The following $\Set$-based structures are finite limit definable:
simplicial set, more generally $\D$-set, for any diagram $\D$; monoid,
group, ring$\dots$ more generally, any species of finitary equational
universal algebras (possibly many-sorted); (small) category, groupoid,
2-category, double category$\dots$, diagram of a fixed shape within the
category of any of the aforementioned structures.  The corresponding
notion of morphism of structures is always (strict) homomorphism,
preserving all operations on the nose.
\end{example}

\begin{remark}
As pointed out by the referee, both \ref{flimfact}(iii) and (iv) are
special cases of the classical fact that given cartesian structures $S_1$
and $S_2$, there exists a ``tensor product'' cartesian theory $S_1\tens
S_2$ such that models of $S_1$ in the category of $S_2$-structures, as well
as models of $S_2$ in the category of $S_1$-structures, are canonically
isomorphic to $S_1\tens S_2$-structures.
\end{remark}

The prototype of (iv) is the observation that a natural transformation
between two diagrams, say of shape $\D$, can itself be thought of as a
diagram, indexed by $\{\bullet\ra\star\}\times\D$.  It reduces the task of
understanding how to specify well-behaved classes of morphisms between
finite limit structures --- notably, the weak equivalences and
(co)fibrations --- to specifying subcategories of models (of a different
theory).  Now recall property (3) of the motivating examples, that weak
equivalences are preserved by inverse images, and the solution set
condition \bd{c3} of Theorem~\ref{jeff}.  These are characteristic of
notions that can be defined in terms of finite limits and arbitrary
colimits, using a set's worth of data.  Since the early 70's, categorical
logicians have provided several equivalent formal analyses of such notions,
a convenient one being

\subsection{Coherent logic} This fragment is richer than cartesian
logic.  The language can contain (besides typed variables and function
symbols, as usual) relation symbols as well as the logical connectives
$\es,\vagy,\implies,\forall,\exists,\top,\bot$.  However, only sentences of
the following form are allowed:
\[             \forall\vec{x}\bigl(\Phi\implies\Psi\bigr)          \]
where $\vec{x}$ stands for a string of variables, and $\Phi$ and $\Psi$
must be built up via the use of $\exists$, finite conjunctions and
arbitrary (i.e.\ possibly infinite) disjunctions from ``atomic formulas'',
which are meaningful combinations of constants, variables and function and
relation symbols.

\begin{example}   \label{cohex}
Within each category of type $\Mod_S(\Set)$ below, $S$ a finite limit
structure, the notions listed are coherently definable in the language of
$S$:
\begin{itemize}
  \item[$\bullet$] $\SSet$: fibrant simplicial sets (i.e.\ those satisfying
the Kan extension condition), acyclic simplicial sets
  \item[$\bullet$] $\Mor(\SSet)$: monomorphisms, Kan fibrations, 
``topological'' weak equivalences
  \item[$\bullet$] Groups: divisible groups, torsion groups
  \item[$\bullet$] Rings: local rings, fields, fields of characteristic
$p>0$, separably closed fields, algebraically closed fields
  \item[$\bullet$] $\Mor$(Rings): flat homomorphisms
  \item[$\bullet$] $\D_+$-diagrams in $\Mod_S$, where $\D_+$ is the cocone
on a diagram $\D$ (i.e.\ is the result of formally adding a disjoint
terminal object to the small category $\D$): being an initial cocone, i.e.\
a colimit on a functor $\D\ra\Mod_S(\Set)$.
  \item[$\bullet$] $\D_-$-diagrams in $\Mod_S$, where $\D_-$ is the cone on
a finite diagram $\D$ (i.e.\ is the result of formally adding a disjoint
initial object to a category $\D$ with finitely many arrows): being a
terminal cone, i.e.\ a limit on a functor $\D\ra\Mod_S(\Set)$.
\end{itemize}
\end{example}

The reader wishing to see coherent logic in action (in addition to an
introduction to it) may enjoy reading Wraith~\cite{wraith77} that
constructs the \'etale topos of a scheme without the intervening use of a
topology on the category of schemes, or any subcategory thereof.  (The
general \'etale topos is glued together from the affines, and the \'etale
topos of a ring is constructed as a classifying topos.)  Johnstone's
related \cite{john77} shows how to obtain various spectral (i.e.\ sheaf)  
representations of rings directly by looking at the axioms and also treats
Coste's logic.

Much as cartesian logic matches finite limit sketches, the language of
coherent logic is interdefinable with colimit-and-finite limit sketches;
see Makkai--Par\'e~\cite{makkpare}, Borceux~\cite{borceux}~vol.III.

CAVEAT.  The terms \emph{geometric} and \emph{coherent} logic are sometimes
used interchangeably.  When not, they mean almost the same thing: one of
them refers to a logic permitting infinitary disjunctions (arbitrary
colimits) while the other to the same logic with only finitary disjunctions
(resp.\ finite colimits).  In this article, the term ``geometric logic'' is
not used and ``coherent logic'' allows arbitrary infinite disjunctions.

The next proposition is a ``reader's digest'' of facts of coherent logic,
tailored for the needs of the present paper.

\begin{prop}   \label{cohefact}
Let $S$ be a structure defined in terms of finite limits, $\E$ a topos,
$\A$ a set of coherent sentences in the language of $S$.  If 
$X\in\Mod_S(\E)$ satisfies these sentences (``axioms''), one writes
$X\models\A$.  Let $\Mod_{S,\A}(\E)$ be the full subcategory of
$\Mod_S(\E)$ with objects those $X$ that satisfy $\A$.
\begin{itemize}
   \item[(i)] Let $\E\llra{f}\F$ be a topos morphism,
$X\in\Mod_S(\F)$.  If $X\models\A$ then $f^*(X)\models\A$.
   \item[(ii)] $\Mod_{S,\A}(\E)$ is closed in $\Mod_S(\E)$ under
filtered colimits.
   \item[(iii)] $\Mod_{S,\A}(\E)$ is an accessibly embedded, accessible
full subcategory of $\Mod_S(\E)$.
   \item[(iv)] Let $\D$ be a small category.  For any $d\in\ob\D$, let
$\ev_d$ denote the ``evaluation at $d$'' functor $\E^\D\ra\E$, which is the
inverse image of a topos morphism.  For $X\in\Mod_S(\E^\D)$, $X\models\A$
iff $\ev_d(X)\models\A$ for all objects $d$ of $\D$.
\end{itemize}
\end{prop}

Example~\ref{cohex} may suggest that one can transplant a Quillen model
category from $\Set$ to an arbitrary topos by borrowing the ``formulaic''
definition of weak equivalences, cofibrations and fibrations.  This is too
much to ask, as pointed out by Jardine~\cite{jard87}: any
Eilenberg--MacLane sheaf is locally fibrant (being a sheaf of simplicial
abelian groups); if it were globally fibrant as well, sheaf cohomology
would be trivial.  Experience has shown that the ``local'' (i.e.\
stalkwise, i.e.\ ``formulaic'' in the sense of the internal logic) weak
equivalences are the right ones.  Any two of the three classes defining a
Quillen model category determine the third.  This suggests that one borrows
the logical description of weak equivalences and one of the classes of
cofibrations and fibrations.  The second possibility, in general, leads to
a local homotopy theory, or ``homotopy theory of fibrant objects'', in the
sense of K.~Brown~\cite{kenbrown}; see also Jardine~\cite{jard86}.  Note
that Quillen's axioms are self-dual; the asymmetry is inherent in the
set-theoretic aspect of sheaf theory (cf.\ cofibrant generation etc.).

The first choice leads to a Quillen model category and is the subject of
this paper.  Note that Jeff Smith's theorem prefigures the non-uniqueness
of the choice of cofibrations: provided one deals with an accessible class
of weak equivalences, which therefore satisfy the solution set condition
everywhere, one may replace the $I$ of \ref{jeff} by any set
$I'\subset\cof(I)$ as long as condition \bd{c1} is not violated.  For this
reason, the main theorem below is split into two parts; the first one
produces cofibrations less uniformly in $\TOPOI$, but under less stringent
conditions.

For the rest of this section, let $S$ stand for the definition of a
structure in terms of finite limits.  For brevity, write $\tS(\E)$ for
$\Mod_S(\E)$, $\E$ a topos.  Let $W$ resp.\ $C$ be two sets of axioms of
coherent logic in the language of morphisms of models of $S$. Write
$\tW(\E):=\{f\in\mor\tS(\E)\;|\;f\models W\}$ and
$\tC(\E):=\{f\in\mor\tS(\E)\;|\;f\models C\}$.

\begin{thm} \label{auto} Consider the hypotheses
\item[(i)] $\tS(\Set)$ with weak equivalences $\tW(\Set)$ and cofibrations
$\tC(\Set)$ is a Quillen model category; and $\tC(\Set)$ is of the form
$\cof(I)$, $I$ a set.
\item[(ii)] For every topos $\E$, $\tC(\E)=\cof(I_\E)$ for some set $I_\E$
of maps.
\item[(iii)] One of the following: $(\dag)$ $W$ and $C$ belong to a 
fragment of coherent logic that has enough models in $\Set$ (for example,
axioms of finite length; the countable fragment of coherent logic;
essentially algebraic theories) or $(\ddag)$ $\tS(\E)$ with weak
equivalences $\tW(\E)$ and cofibrations $\tC(\E)$ is a Quillen model
category for every topos $\E$ of the form $\sh(\B)$, $\B$ being a complete
Boolean algebra with its canonical topology.

\item[\bd{(1)}] (i) implies that for every topos $\E$ with enough points,
$\tS(\E)$ with weak equivalences $\tW(\E)$ and cofibrations a certain
subclass of $\tC(\E)$ is a cofibrantly generated Quillen model category.  
In the presence of (iii), the conclusion extends to every topos. 

\item[\bd{(2)}] (ii) implies that for every topos $\E$ with enough points,
the cofibrations can be chosen to be $\tC(\E)$, the rest being as in
\textup{\bd{(1)}}.  In the presence of (iii), the conclusion extends to
every topos.
\end{thm}

\begin{remark}
It is immediate that in case \bd{(2)}, a topos morphism 
$\E\three{f^*}{\leri}{f_*}\F$ induces a Quillen pair $\tS(\E)\leri\tS(\F)$.
(\ref{flimfact}(ii) and \ref{cohefact}(i) show that the left adjoint
preserves \emph{all} weak equivalences and cofibrations.)  One thus has a
``coherently definable'' or ``sheafifiable'' homotopy theory and a functor
$\TOPOI\ra\HOMODEL$.  See remark~\ref{topshape} as to functoriality under
\bd{(1)}.
\end{remark}

\begin{remark}
The choice of cofibrations is not unique, even if all assumptions are
satisfied.  There may be more than one functorial cofibration class as
well.
\end{remark}

\begin{remark}
In all examples I am aware of, if (i) holds, so does (ii).  In fact, (ii)
is likely to be a formal consequence of its being true for the case
$\Set=\E$.  Alas, as things stand now, (ii) has to be checked by hand,
unless \ref{mono} applies.
\end{remark}

\begin{remark}
I know a single example where one doesn't seem to get by with countably
many defining axioms of countable length, or with some set of axioms of
finite length, for $W$ and $C$ --- cases when (iii)$(\dag)$ applies --- and
that is homological localization of simplicial sets (simplicial objects in
a topos).  That example is much better seen as a direct consequence of Jeff
Smith's theorem, together with \ref{mono} and a bit of accessibility.
\end{remark}

The statement of \ref{auto} begs its proof: observe that the hypotheses of
Jeff Smith's theorem hold for $\tS(\Set)$; transfer them to $\tS(\E)$ via
logical methods; forward applications of Thm.~\ref{jeff} yield the
conclusions.  To begin with, one has unconditionally

\begin{lemma}   \label{condi1}
For any topos $\E$, $\tS(\E)$ is a locally presentable category.  $\tW(\E)$
is closed under retracts and satisfies the solution set condition in
$\Mor(\tS(\E))$.
\end{lemma}

\begin{proof}
By \ref{flimfact}(i), \ref{cohefact}(iii), \ref{retrac} and
\ref{accsolset}.
\end{proof}

The next two lemmas state facts in tandem, preceded by the hypotheses they
need. 

\begin{lemma}  \label{condi2}
\item[(i)] For every topos $\E$ with enough points: $\tW(\E)$ has the
2-of-3 property; $\tC(\E)$ is closed under composition; 
$\tC(\E)\cap\tW(\E)$ is closed under pushout.
\item[(i,iii)] Same conclusions for \emph{every} topos $\E$.
\end{lemma}

\begin{proof}
Each of these properties has the following form: a set of coherent
sentences (in the language of the appropriate diagram of $S$-structures)
implies certain coherent sentences.  Topoi with enough points inherit the
truth of such statements from $\Set$.  As to $(\dag)$, the definition of
having enough models in $\Set$ is that such conclusions extend to an
arbitrary topos.  It is a theorem of Makkai--Reyes that countable coherent
logic has enough models in $\Set$, of Deligne-Joyal that finitary coherent
sentences do.  The case of universal Horn logic / essentially algebraic
theories is classical (see Makkai--Reyes~\cite{makkreyes}).  The Boolean
case $(\ddag)$ is Barr's theorem.
\end{proof}

\begin{lemma}   \label{condi3}
\item[(i)] For every topos $\E$ with enough points: $\tW(\E)$ and $\tC(\E)$
are closed under transfinite composition.
\item[(i,iii)] Same conclusions for \emph{every} topos $\E$.
\end{lemma}

\begin{proof}
By Prop.~\ref{cohefact}, $\tW(\E)$ and $\tC(\E)$ are closed under filtered
colimits in the category of morphisms of $\tS(\E)$.  If they are closed
under composition, they must be closed under transfinite composition.  (Use
transfinite induction.)
\end{proof}

\begin{lemma}  \label{condi4}
(i) implies: for any topos $\E$, $\inj(\tC(\E))\subseteq\tW(\E)$.
\end{lemma}

\begin{proof}
This holds for $\E=\Set$ by assumption, since acyclic fibrations are weak
equivalences.  Extend the conclusion from $\Set$ to presheaf topoi 
$\pre(\D)$.  Let $d$ be any object of the small category $\D$, and consider
the adjunction
\[
    \tS(\Set)^{\D^\op} \three{L}{\leri}{\ev_d} \tS(\Set)
\]
where $\ev_d$ is the ``evaluation at $d$'' functor and $L$ its left adjoint
(the left Kan extension).  Let $f\in\tC(\Set)$.  $L(f)$ is a 
$\D^\op$-diagram in $\tS(\Set)$ that is, at every object of $\D$, a copower
of $f$.  Recall that $\tS(\Set)^{\D^\op}$ is the ``same'' as 
$\tS(\pre(\D))$.  Since any copower of $f$ is a cofibration in $\tS(\Set)$,
and since coherent axioms are evaluated ``objectwise'' in functor
categories (cf.~\ref{cohefact}(iv)), $L(f)$ will belong to the class
$\tC(\pre(\D))$.  Suppose $g\in\inj(\tC(\pre(\D))$.  A fortiori 
$g\in\inj(L(\tC(\Set)))$.  By adjunction, $\ev_d(g)\in\inj(\tC(\Set))$, so
$\ev_d(g)\in\tW(\Set)$.  But since the class $\tW(\pre(\D))$ is defined
again by coherent axioms, $\ev_d(g)\in\tW(\Set)$, for every object $d$ of
$\D$, implies $g\in\tW(\pre(\D))$.

Consider now an arbitrary topos $\E$.  Choose a site $(\D,J)$ of definition
for $\E$, and consider the inclusion $\E\three{\ell}{\leri}{i}\pre(\D)$,
where $\ell$ is sheafification.  It induces an adjunction 
$\tS(\E)\leri\tS(\pre(\D))$.  (Out of laziness, retain the same
letters to denote these adjoints.)  Take any $f\in\inj(\tC(\E))$.  Since
sheafification (being an inverse image part of a topos morphism) preserves
the coherently defined class of cofibrations, by adjointness one has
$i(f)\in\inj(\tC(\pre(\D)))$ whence $i(f)\in\tW(\pre(\D))$.  Since
sheafification must preserve weak equivalences too, $f\iso\ell
i(f)\in\tW(\E)$.
\end{proof}

Part \bd{(2)} of \ref{auto} is now simply Jeff Smith's theorem~\ref{jeff}
applied to the data $\tS(\E)$, $\tW(\E)$, and the $I_\E$ of assumption
(ii), using Prop.~\ref{flimfact}(i) and lemmas \ref{condi1} through
\ref{condi4}.

Part \bd{(1)} follows by a closer look at the argument in \ref{condi4}.
Assume (i).  By Jeff Smith's theorem, $\tS(\Set)$ is a cofibrantly
generated Quillen model category.  By a theorem of
Dwyer--Hirschhorn--Kan~\cite{hirsch}, for any cofibrantly generated model
category $\M$ and small category $\D$, $\M^\D$ possesses a model structure
with the weak equivalences and fibrations being $\D$-objectwise.  
Cofibrations in $\M^\D$ can be described as follows.  Let $\D_\delta$ be
the diagram $\D$ ``made discrete'', that is, $\D_\delta$ is the set of
objects of $\D$ and their identities.  There is an adjunction
$\M^\D\leri\M^{\D_\delta}$ where the right adjoint is the forgetful functor
(i.e.\ precomposition with $\D_\delta\ra\D$) and the left adjoint its left
Kan extension $L_K$.  $M^{\D_\delta}$ is (tautologously) a cofibrantly
generated model category whose set of generating cofibrations,
$I_{\D_\delta}$, is the set of tuples that are coordinatewise generating
cofibrations in $\M$.  With the above choice of weak equivalences and
fibrations, cofibrations in $\M^\D$ are $\cof(L_K(I_{\D_\delta}))$.

By virtue of \ref{flimfact}(iii) and \ref{cohefact}(iv), the above
prescription gives a cofibrantly generated Quillen model structure on
$\tS(\pre(\D))$ with weak equivalences the $\tW(\pre(\D))$ and cofibrations
a --- for non-trivial $\D$, proper --- subclass of $\tC(\pre(\D))$.  (This
follows also by a direct application of Jeff Smith's theorem.)  Now for an
arbitrary topos $\E$, choose a site, i.e.\ topos inclusion
$\E\three{\ell}{\leri}{i}\pre(\D)$.  Write, for brevity, $K$ for the set of
generating cofibrations in $\tS(\pre(\D))$.  Apply Jeff Smith's theorem to
the data $\tS(\E)$, $\tW(\E)$, $\ell(K)$.  Since $\tS(\E)$ is a full,
reflective subcategory of $\tS(\pre(\D))$, $\cof(\ell(K))=\ell(\cof(K))$.  
That is, cofibrations in $\tS(\E)$ are sheafifications of the cofibrations
in $\tS(\pre(\D))$.  So $\cof(\ell(K))\subseteq\tC(\E)$, which lets one
deduce that a pushout of an acyclic cofibration is acyclic from
lemma~\ref{condi2}.  Use \ref{condi1} and \ref{condi3} for \bd{c0}, \bd{c3}
and the other part of \bd{c2}.  The adjunction argument used in the second
part of the proof of \ref{condi4} establishes \bd{c1}.  This completes the
proof of \ref{auto}.  \qed

\begin{example}   \label{bkmono}
Here's how to construct a model category on simplicial objects in a topos
with Joyal's weak equivalences, whose class of fibrations is intermediate
between the local fibrations (or fibrations in the internal sense; when the
topos has enough points, these are the maps that are stalkwise fibrations)
and the global ones, i.e.\ maps with the right lifting property w.r.t.\
every acyclic mono.  Say $\E=\sh(\D,J)$.  Consider the Quillen model
structure on $\pre(\D)^\si=\SSet^{\D^\op}$ created by the forgetful functor
to $\SSet^{\D_\delta}$.  (It is the Bousfield--Kan model structure, see
\cite{bouskan} XI.8.)  Its class of cofibrations $\bk$ is a proper subclass
of the monomorphisms in general.  Letting 
$\E\three{\ell}{\leri}{i}\pre(\D)$ be the canonical adjunction, define 
cofibrations in $\E^\si$ to be the class $\ell(\bk)$, weak equivalences
$\tW$ as in \ref{joy} and fibrations to be $\inj(\ell(\bk)\cap\tW)$.
Proof: choose a set $I$ such that $\cof(I)=\bk$ in $\pre(\D)^\si$.  Then
$\ell(\bk)=\cof(\ell(I))$ in $\E^\si$.  Apply Jeff Smith's theorem.

If the topology $J$ is ``too strong'', $\ell(\bk)$ may include all
monomorphisms in $\sh(\D,J)$.  (For example, should it happen that
$\sh(\D,J)$ is equivalent to $\Set$, note that $\cof(I)$ is all injections
as soon as $I$ contains a single non-trivial injection.)  There seems to be
no a priori reason for this to happen for all non-trivial topologies,
however.
\end{example}

\begin{remark}  \label{topshape}
This formal argument extends to the upcoming examples, and even to the ones
in \cite{right}.  Despite the unsettling vagueness of the ``correct''
cofibration class, any two choices have a common super-class, thus giving
rise to Quillen-equivalent model categories \cite{cofib}.  In fact, the
Quillen equivalence type of any presentable model category is fixed by the
category itself and its subcategory of weak equivalences.  This suggests
that the proper target of $\TOPOI\ra\HOMODEL$ may be the category of
Quillen model categories and Quillen pairs \emph{modulo} Quillen
equivalence.  Note that existing work of Dwyer--Kan and Rezk also suggest
that ``a homotopy theory'' (as such) is determined by the category of
models and its subcategory of weak equivalences.
\end{remark}

\section{Examples} As soon as there is one example of a sheafifiable
homotopy model category --- in the sense of part \bd{(2)} of \ref{auto} ---
there are infinitely many, by the following cheap observation: since
$\tS(\E)^\D$ is equivalent to $\tS(\E^\D)$, and $\E^\D$ is a topos if $\E$
is, there exists a sheafifiable homotopy theory of $\D$-diagrams of
$S$-structures with the weak equivalences and \emph{cofibrations} being
$\D$-objectwise.  This is reminiscent of (and contains as a special case)  
the model structure on simplicial diagrams that Heller~\cite{heller}
denotes \emph{right}: let the cofibrations be all the monos.  It is in
general distinct from (but Quillen-equivalent to) the one existing on
diagrams over any cofibrantly generated model category (reducing to the
Bousfield--Kan structure on simplicial diagrams that Heller named
\emph{left}).  We list six non-trivial cases of the main theorem here.  
Each satisfies all three hypotheses of \ref{auto}; (ii) holds either
because Prop.~\ref{mono} applies directly, or because the structure arises
via an adjunction from one where it applies.  This is examined in more
detail in \cite{right}.  We also digress into the unbounded derived
category.

\begin{example} (simplicial sets) \\    \label{sset}
The observation that weak equivalences of simplicial sets are definable in
terms of finite limits and countable colimits goes back (at least) to
Illusie~\cite{illusie}.  It is not necessary to construct the homotopy
group objects internally, however.  To begin with, consider a map
$\Xb\llra{f}\Yb$ between Kan complexes in $\SSet$.  $f$ is a weak
equivalence iff it induces a bijection between embedded ``singular
spheres'', that is, simplices all of whose lower-dimensional faces are one
and the same 0-simplex and its degeneracies.  That the map is onto amounts
to: \emph{for every 0-simplex $y_0$ in $\Yb$ and $n$-simplex $y_n$ all of
whose faces are $y_0$ (and its degeneracies), there exists an $n$-simplex
$x_n$ all of whose faces are some $x_0$ (and its degeneracies) such that
$f(x_n)$ is in the same based homotopy class of singular simplices as
$y_n$}, where for two $n$-simplices to be in the same homotopy class
amounts to the existence of an $n+1$-simplex with suitable face matching
conditions.  The sentence  in italics translates verbatim into a coherent
axiom in the language of morphisms of simplicial objects.  That $f$
induces an injection is similar.  To deal with an arbitrary simplicial
map, one has to use some fibrant replacement functor $\SSet\ra\SSet$
definable in coherent logic.  Kan's \cite{kan57} $\exfty$ is such.  The
functor $\Ex$ --- together with the natural transformation $\Id\ra\Ex$ ---
is defined in terms of finite limits, and $\exfty$ is the colimit along
the (countable) chain of iterations of $\Id\ra\Ex$.
\end{example}

\begin{example} (homological localizations of simplicial sets) \\
Let $h_*$ be a homology theory on $\SSet$.  We need to choose a
representation of $h_*$ in the sense of G.\ Whitehead~\cite{whit62},
\begin{equation}
  h_n(X) = \colim \pi_n(X_+ \es E_i)      \label{holoc}
\end{equation}
where $E_i$ is a ``naive spectrum'', ie.\ sequence of pointed simplicial
sets and connecting maps from the suspension of $E_i$ to $E_{i+1}$.  One
may describe by coherent axioms (or, colimit and finite limit constructions
of simplicial sets) the following, in turn: adding a disjoint basepoint to
$X$; smashing two pointed simplicial sets; suspending a simplicial set (one
ought to use Kan's model here); computing $\pi_n(-)$ of a simplicial set
(after fibrant replacement); the $\colim_i$ of homotopy groups (here one
has to use that the connecting maps are definable too); finally, that for a
simplicial map $X\llra{f}Y$, $h_n(f)$ is an isomorphism.

Theorem~\ref{jeff} and Prop.~\ref{mono} can now be applied; this is one of
the (rare!)\ cases when \bd{c2} can be checked directly.  Note that by 
using the machinery of accessible functors and classes, Bousfield's
original proof in \cite{bous75} can also be made to work in any topos
\cite{tezis}.  By working directly in the site, Goerss and Jardine
\cite{goerssjard} establish homological localization for simplicial
sheaves as well.
\end{example}

\begin{remark}
The calculus of accessibility highlights one aspect of the
homology--cohomology dichotomy.  A homology theory is a functor
$\SSet\llra{h_*}\A$ where $\A$ is, typically, a locally presentable
category.  ($\A$ will be so as soon as it is a cocomplete $AB5$ abelian
category with generator --- such as, graded abelian groups.  See
Prop.~\ref{grolo}.)  By Milnor's axiom, $h_*$ preserves filtered colimits.
In fact, as soon as $h_*$ preserves $\kappa$-filtered colimits for some
$\kappa$, it is an accessible functor, and $h_*^{-1}(\isom)$ will be
accessible in $\Mor(\SSet)$, yielding the solution set condition.  By
contrast, a cohomology theory is a functor $\SSet\llra{h^*}\A^\op$.  If the
opposite category of a locally presentable category is accessible, they
must be (equivalent to) small posets; see \AdRos~1.64.  So $\A^\op$ is not
accessible; the accessibility of $h_*^{-1}(\isom)$ in $\SSet$ is then an
open question.

From the set-theoretical point of view, all that was used of $\SSet$ is
that it is a locally presentable category.  Indeed, J.\ Smith's theorem 
implies that any presentable model category allows localizations w.r.t.\
any accessible homology theory.
\end{remark}

\begin{remark}
This is a good place to discuss (in)dependence of the choices made.  After
all, the class of weak equivalences in $\SSet$ admits a canonical 
definition: the maps whose topological realizations are weak homotopy
equivalences.  In Ex.~\ref{sset}, \emph{some} combinatorial 
characterization of this had to be found.  Similarly, in Ex.~\ref{holoc},
\emph{a} representing spectrum had to be chosen.  Dependence on these
choices is not a frivolous issue; there are examples (albeit very
artificial, from the point of view of homotopy theory) of coherent axioms
that cannot be satisfied in $\Set$ (so, from the point of view of sets,
parametrize the empty collection) but do have models in other topoi.  On
the positive side, one has
\end{remark}

\begin{prop}
Let $A_1$, $A_2$ be two coherent axiom systems in the same (first-order,
many-sorted) language $\L$.  Suppose that in the category of sets, the same
structures satisfy them.
\begin{itemize}
\item[(i)] Then the same is true in any topos with enough points.
\item[(ii)] Suppose that there is a ``reason'' for this coincidence, more
precisely, a formal deduction via coherent logic of $A_1$ from $A_2$, and
vice versa.  Then the same class of structures satisfies them in any topos.
\item[(iii)] Suppose that both $A_1$ and $A_2$ belong to a fragment of
coherent logic that has enough models in $Set$.  (For example, $A_i$,
$i=1,2$, may be finitary, or employ countable disjunctions but have only
countably many axioms, or be an essentially algebraic theory.)  Then again,
the same structures will satisfy them in any topos.
\end{itemize}
\end{prop}

\begin{proof}
(The reader looking for relevance to abstract homotopy theory may wish to
skip to the examples following the proof.)  (i) is trivial, given that
points of a topos preserve and any conservative collection of points
reflect the truth of coherent axioms.  (ii) is the mere fact that coherent
deductions remain valid in any topos.  (iii) is a ``degenerate case'' of
the Conceptual Completeness theorem of Joyal--Makkai--Reyes~\cite{makkreyes}.
$A_1$ and $A_2$, by assumption, have enough models and have the same models
in $\Set$.  Hence they have the same class of coherent consequences.  This
means that their associated syntactic sites $\Def(A_1)$, $\Def(A_2)$ (see
Makkai--Reyes~\cite{makkreyes} or MacLane--Moerdijk~\cite{m&m}) are
equivalent.  (In fact, $\Def(A_1)$ and $\Def(A_2)$ can be constructed to be
the \emph{same} --- not only isomorphic, but identical --- since $A_1$ and
$A_2$ were assumed to share the same language $\L$.)  But this implies they
have the same models in any topos.
\end{proof}

(i) already covers all but a few types of topoi that tend to arise in
algebraic geometry and topology.  (ii) applies, for example, to the choice
of representing spectra.  Let $E_1$ and $E_2$ be weakly homotopy equivalent
spectra, and let $A_1$ say ``it (i.e.\ a variable spectrum) has the same
weak homotopy type as $E_1$''; analogously for $A_2$.  Should these two
classes coincide, there is a witness for that, namely a \emph{map of
spectra $E_1\llra{f}E_2$ inducing isomorphisms on the stable homotopy
groups}.  The statement in italics is equivalent to a set of sentences of
coherent logic.  This allows the formal demonstration of the equivalence of
``the spectrum $X$ has the same weak homotopy type as $E_1$'' with ``the
spectrum $X$ has the same weak homotopy type as $E_2$'' valid in any topos. 
(iii) implies, for example, that any ``purely combinatorial'' definition of
a simplicial map being a weak equivalence can be chosen as long as it has
the intended meaning for simplicial \emph{sets}.  The technical sense of
``purely combinatorial'' is: statements in the countable fragment of
coherent logic, with signature $\Delta$ (and its extension by coproducts
and coequalizers of equivalence relations).  Intuitively, the definition
must operate directly with simplices and the face and degeneracy maps (no
mapping spaces, fundamental groupoids etc.\ can be used unless these had
been so defined beforehand).  The axioms must be ``if--then'' type,
demanding the existence of simplices, or sequences of simplices, satisfying
finitary face-matching conditions whenever other conditions of this type
are met; but one cannot employ uncountably many conditions nor uncountable
strings of simplices.

\begin{example} (small categories and categorical equivalences) \\
Motivating example~\ref{sstack} is a case of Theorem~\ref{auto} where not
all the monomorphisms are cofibrations.  I am indebted to Charles Rezk for
pointing out that this model structure arises from the Reedy model
structure on simplicial spaces (i.e.\ simplicial objects in $\SSet$) via a
right adjoint; this gives a simple proof of condition (ii) of \ref{auto}.
\label{smallcat}
\end{example}

\begin{example} (groupoids and categorical equivalences) \\
Restrict the previous example to groupoids, leaving the definitions of
cofibrations and weak equivalences the same.  (Fibrations are functors
with the right lifting property w.r.t.\ the inclusion 
$\{\bullet\}\inc\{\bullet\ra\star\}$.)  This model structure too was
sheafified ``by hand'' in Joyal--Tierney~\cite{joytier90}.  Via the nerve
functor, the existence of this model category can also be deduced from that
on simplicial objects; this proves \ref{auto}(ii).
\label{grpd}
\end{example}

\begin{example} (cosimplicial simplicial objects, Reedy model structure) \\
For indexing categories $\D$ satisfying a certain combinatorial property,
work of D.\ Kan (see \cite{dhk} or Hovey~\cite{hovey}) shows the existence
of a Quillen model structure on $\M^\D$ with weak equivalences objectwise 
but neither cofibrations nor fibrations so.  For the case $\D=\Delta$ and
$\M=\SSet$, this model structure made its appearance in
Bousfield--Kan~\cite{bouskan} X.4.  It has as cofibrations those
monomorphisms in $\SSet^\Delta$ that induce isomorphisms on the ``maximal
augmentation''.  Here the maximal augmentation of
$X^\bullet\in\ob\SSet^\Delta$ is a subobject of $X^0$, the equalizer of the
coface maps $d^0$, $d^1$.  This is a coherently definable condition.  That
\ref{auto}(ii) is satisfied follows by Kan's methods, or one can use the
following mild extension of Prop.~\ref{mono}: keep hypotheses (i), (ii) and
let $\tC$ be an accessible subclass of the monomorphisms closed under
pushout and transfinite composition and satisfying the following
cancellation condition: if $gf\in\tC$, $f\in\tC$ and $g\in\mono$, then
$g\in\tC$.  Then $\tC=\cof(I)$ for some set $I$.  (Accessibility lets one
find a set of generators; after that, the argument is the same.
See~\cite{tezis}.)
\end{example}

\begin{example} (cyclic sets, ``weak'' model structure) \\
Let $\Lambda^\op$ be Connes' indexing diagram for \emph{cyclic objects},
containing the simplicial indexing diagram $\Delta^\op$ as a subcategory.
Dwyer--Hopkins--Kan~\cite{dhk85} show the existence of a Quillen model
structure on $\pre(\Lambda)$, or ``cyclic sets'', created by the forgetful
functor $\pre(\Lambda)\ra\SSet$.  That is to say, $f\in\mor\pre(\Lambda)$
belongs to $\tW$ iff it is a weak equivalence considered just as a map of
simplicial sets.  Thanks to Ex.~\ref{sset}, this makes the class $\tW$
coherently definable.  Dwyer--Hopkins--Kan define a cyclic map to be a
fibration iff it is a Kan fibration when considered just as a map of
simplicial sets, and give the following combinatorial description of
cofibrations: $X\llra{f}Y\in\pre(\Lambda)$ is a cofibration iff it is an
injection such that for each object \bd{[n]} of $\Lambda^\op$, its group of
automorphisms (which is a cyclic group of order $n+1$) acts freely on the
elements of $Y(\bd{[n]})$ not in the image of $X(\bd{[n]})$.  Taking the
contrapositive, this becomes a coherent condition.
\end{example}

Spali\'nski's \cite{spali} strong model category structure on cyclic sets
is dealt with in \cite{right}.

Another application of Theorem~\ref{jeff} and Prop.~\ref{mono} concerns the
unbounded derived category.  It is valid for Grothendieck abelian
categories, which are automatically locally presentable.  Though this fact
is a simple conjunction of theorems of Gabriel--Popescu and Gabriel--Ulmer,
we write it out in longhand, since the identification is not commonly made
in the literature.

\begin{prop}     \label{grolo}
Let $\A$ be an abelian category that satisfies Grothendieck's axiom
$AB5$ (``directed colimits are exact'', i.e.\ filtered colimits commute
with finite limits). Then the following are equivalent:
\begin{itemize}
\item[(i)] $\A$ is cocomplete and has a generator, i.e.\ object
$G$ such that $\hom(G,-)$ is faithful \footnote{MacLane~\cite{cwm} suggests
the name \emph{separator} for such an object, which is certainly
reasonable, given the plenitude of senses of ``generator'' in different
categorical contexts.}
\item[(ii)] $\A$ is locally presentable.
\end{itemize}
\end{prop}

Historically, an $AB5$ category satisfying (i) has been called a
\emph{Grothendieck (abelian) category}.

\begin{proof}
(ii)$\implies$(i) is trivial: if $\A$ is locally presentable, then it is
(by definition) cocomplete and has a set of dense generators, a fortiori a
set of generators.  Their coproduct will do as a single generator (in the
sense of the proposition), since $\A$ has a zero object.

(i)$\implies$(ii): by the Gabriel-Popescu theorem there exists a ring $R$
and an adjoint pair $\A\three{L}{\leri}{i}\Mod_R$ where $i$ is full and
faithful, and the left adjoint $L$ is exact.  That is, $A$ is equivalent to
a localization of $\Mod_R$.  (Here a \emph{localization} of a category is a
full, reflective, isomorphism-closed subcategory with the reflector
preserving finite limits.)  Such localizations biject with Gabriel
topologies on $R$, i.e.\ collections $\F$ of right ideals of $R$ with
certain closure properties, via associating with $\F$ the full subcategory
of $\Mod_R$ of $\F$-closed modules.  (See e.g.\ Stenstr\"{o}m~\cite{sten}.)
A module $M$ is $\F$-closed iff it is orthogonal to the maps $I\inc R$, 
$I\in\F$, i.e.\ $\hom_R(R,M)\ra\hom_R(I,M)$ is an isomorphism for all
$I\in\F$.  $\Mod_R$ is locally presentable; by the theorem of orthogonal
reflection in locally presentable categories (see e.g.\ \AdRos~Cor.~1.40),
so is its full subcategory of $\F$-closed modules.
\end{proof}

\begin{remark}
The second implication is nontrivial.  $\A$ being locally presentable
entails, for example, that every object of $\A$, in particular, the
generator, has a \emph{rank}, that is, $\hom(G,-)$ commutes with
$\kappa$-filtered colimits for some $\kappa$.  This is not even implicit in
the definition of a generator (as ``separator'').
\end{remark}

\begin{remark}
Prop.~\ref{grolo} is the additive analogue of the case of a topos.  
Giraud's theorem says that if a category $\E$ is cocomplete, has
well-behaved coproducts and quotients of equivalence relations, and a set
of generators, then it is a topos.  It is then equivalent to a
localization of a functor category $\pre(\C)$, where $\C$ is small.  
Localizations of $\pre(\C)$ biject with Grothendieck topologies on $\C$,
via passing to sheaves on a topology; and the sheaf condition can be
phrased as orthogonality w.r.t.\ a certain set of maps in $\pre(\C)$.  
The conclusion that $\E$ is a locally presentable category now follows as
before.
\end{remark}

\begin{prop}   \label{unbounded}
Let $\A$ be a Grothendieck abelian category.  There is a Quillen model
structure on $\Ch_\Z(\A)$ (i.e.\ unbounded chain complexes in $\A$) where
cofibrations are the monomorphisms, and weak equivalences the 
quasi-isomorphisms.
\end{prop}

\begin{proof}
Apply Thm.~\ref{jeff}.  If $\A$ is Grothendieck abelian, so is
$\Ch_\Z(\A)$.  The homology functor $\Ch_\Z(\A)\ra\A^\Z$ (here $\Z$ stands
for just the discrete set) commutes with filtered colimits and $\A^\Z$ is
locally presentable.  The class of isomorphisms is accessible in any
locally presentable category.  Prop.~\ref{inv} gives \bd{c3}, and
Prop.~\ref{mono} grants the generating cofibrations.
\end{proof}

\begin{cor}  \label{exists}
The unbounded derived category of a Grothendieck abelian category exists
(i.e.\ has small hom-sets).
\end{cor}

\begin{remark}
Even when a localization problem $\C[W^{-1}]$ allows a calculus of
fractions, a genuine set-theoretical difficulty remains in establishing
$\C[W^{-1}]$ has small hom-sets (one has to exhibit cofinal small 
subcategories of the large, filtered index categories that arise).  This 
may be the reason that, to the best of this author's knowledge, a proof of
\ref{exists} in this generality has not appeared in print yet.  It is
certainly a folk theorem, however; see e.g.\ Joyal~\cite{joyal}.  A recent
preprint of Tarr\'{\i}o--L\'opez--Salorio~\cite{talosa} gives a
triangulated proof. 

It is quite ironic to observe that while the simple (inductive) injective
replacement arguments break down for unbounded complexes, the proof
strategy of \ref{unbounded} can only apply to a category of complexes that
is cocomplete.  The fibrant replacement functor it produces in abstracto is
quite horrible; it seems that if one wishes to work with explicit
resolutions of unbounded complexes, one should find a judicious sheaf
representation of the category first.
\end{remark}

\begin{remark} 
In his Toh\^oku classic, Grothendieck was able to demonstrate remarkable
(ultimately, set-theoretical) features of what are now called Grothendieck
abelian categories, notably the existence of enough injectives (via
showing, essentially, that the class of monomorphisms is generated by a
set).  With hindsight, when working with Grothendieck abelian categories,
one encounters a class of examples with much stronger set-theoretic bounds
than is apparent from the concept of generator (``separator'').  When
limits and colimits are not as well-behaved as for abelian categories and
topoi, it seems necessary to posit the extra set theoretic control 
explicitly.  Conversely, perhaps it is not frivolous to say that 
presentable Quillen model categories are ``convenient categories to do
homotopical algebra in'', and to view them as non-abelian counterparts of
Grothendieck abelian categories.
\end{remark}

As a closing remark, note that \cite{right} lists about a dozen Quillen
model categories whose weak equivalences are coherently definable and whose
cofibrations are of the form $\cof(I)$, where $I$ is the image under a
suitable left adjoint of a coherently defined class of maps.  (They thus
generalize examples \ref{smallcat} and \ref{grpd}.)  Their qualitative
features --- sheafifiability, existence of non-canonical as well as
functorial choices for cofibrations, and the essential uniqueness of
cofibrations --- are so similar to those above that it did not seem
worthwhile to burden this paper with the extra logical machinery needed to
``automate'' their construction.

\end{document}